\documentclass{birkjour}
 \newtheorem{thm}{Theorem}[section]

 \theoremstyle{definition}
 
 \theoremstyle{remark}

 \numberwithin{equation}{section}

\usepackage{algorithmic} 
\usepackage[Algorithm, subsection]{algorithm} 
\usepackage{graphicx} 
\usepackage{subfig} 

\begin{document}

\def\ds{\,\mathrm ds}
\def\dt{\,\mathrm dt}
\def\dx{\,\mathrm dx}
\def\b#1{\boldsymbol{#1}}
\def\Epsilon{\mathcal{E}}
\def\R{\mathbb{R}}

\newcommand{\dd}{\mathrm{d}}
\newcommand{\Q}{\Omega_T}
\newcommand{\y}{\boldsymbol c} 
\renewcommand{\u}{\boldsymbol f} 
\newcommand{\p}{\boldsymbol p} 
\newcommand{\bdelta}{\boldsymbol \delta} 
           
\def\N{\mathbb{N}}
\def\der{\mathrm{D}}

%
%
%
%
%
%
%
%
%

\title[]
 {Optimal control of Allen-Cahn systems}

\author[]{Luise Blank, M. Hassan Farshbaf-Shaker, Claudia Hecht, Josef Michl, Christoph Rupprecht}



\date{\today}


\begin{abstract}
Optimization problems governed by Allen-Cahn systems including elastic effects are formulated and first-order necessary optimality conditions are presented. Smooth as well as obstacle potentials are considered, where the latter leads to an MPEC. Numerically, for smooth potential the problem is solved efficiently by the Trust-Region-Newton-Steihaug-cg method. In case of an obstacle potential first numerical results are presented.
\end{abstract}

\maketitle
\section{Introduction and problem formulation}\label{Introduction}
Optimization problems with interfaces and free boundaries frequently appear in materials science, fluid dynamics and biology (see i.e. \cite{Bossavit} and the references therein). In this paper we concentrate on a phase field approach, more precisely on a multi-component Allen-Cahn model, to desribe the dynamics of the interface. This allows complex topological
changes. The possibly sharp interface between the phases is replaced
by a thin transitional layer of width $\mathcal O(\varepsilon)$ where $\varepsilon>0$ is a small parameter, and the $N$ different phases are described by a phase field variable $\b c=(c_1,\ldots,c_N)^T$, where $c_i$ denotes the fraction of the $i$-th material. The underlying non-convex interfacial energy is based on the generalized Ginzburg-Landau energy, see \cite{G2000},
\begin{align}
E({\boldsymbol{c},\boldsymbol{u}}):=\int_{\Omega}\left\{\frac{\varepsilon}{2}|\nabla{\bf c}|^{2}+\frac1\varepsilon\Psi({\bf c})+W(\bf{c},{\mathcal E}(\boldsymbol u))\right\}\dx,\label{Ginzburg-Landau-Energy}
\end{align}
where $\Omega\subset\mathbb R^{d}$, $1\leq d\leq 3$, is a bounded domain with either convex or $C^{1,1}$-boundary. Moreover, ${\boldsymbol u}$ is the displacement field mapping into $\R^d$ and $\Psi$ is the bulk potential. In general the potential $\Psi$ is assumed to have global minima at the pure phases and in physical situations there are many choices possible, see \cite{Blowey-Elliott}. Here we consider two different cases: a smooth double-well potential in Sections \ref{smooth} and \ref{Christoph}, and a nonsmooth obstacle potential in Sections \ref{Obstacle} and \ref{Luise}. The latter ensures in particular that the pure phases correspond exactly to $c_i=1$, whereas in the smooth case those are given by $c_i\approx 1$. The term $W(\bf{c},{\mathcal E}(\b u))$ in $\eqref{Ginzburg-Landau-Energy}$ is the elastic free energy density. Since in phase separation processes of alloys the deformations are typically small we choose a theory based on the linearized strain tensor (see \cite{ciarlet1988three}) given by ${\mathcal E}:={\mathcal E}(\boldsymbol u)=\frac12(\nabla\b u+\nabla\b u^T)$ and
\begin{align}
W(\b c,{\mathcal E})=\frac{1}{2}({\mathcal E}-{\mathcal E}^*(\b c)):{\mathcal C}({\mathcal E}-{\mathcal E}^*(\b c)).\label{mechanical-energy}
\end{align}
Here ${\mathcal C}$ is the symmetric, positive definite, possibly anisotropic elasticity tensor mapping from symmetric tensors in $\mathbb R^{d\times d}$ into itself. The quantity ${\mathcal E}^*(\boldsymbol c)$ is the eigenstrain at concentration $\boldsymbol c$ and following Vegard's law we choose ${\mathcal E}^*(\boldsymbol c)=\sum_{i=1}^N c_{i}{\mathcal E}^*(\boldsymbol e_{i})$, where ${\mathcal E}^*(\boldsymbol e_{i})$ is the value of the strain tensor when the material consists only of component $i$ and is unstressed. Here $(\b e_i)_{i=1}^N$ denote the standard coordinate vectors in $\R^N$. The dynamics of the interface motion can be modelled by the steepest descent of (\ref{Ginzburg-Landau-Energy}) with respect to the $L^2$-norm, see \cite{blesgen,G_habil}. The mechanical equilibrium is obtained on a much faster time scale and therefore we assume quasi-static equilibrium for the mechanical variable $\b u$. For a smooth potential $\Psi$ this results after suitable rescaling of time in the following elastic Allen-Cahn equation
\begin{align}\label{elastic-Allen}
\left( \begin {array} {c}
        \varepsilon\partial_{t}\boldsymbol c \\
        \boldsymbol 0\\         
    \end {array} 
\right)
=
\left( \begin {array} {c}
        \varepsilon\Delta\b c-\frac{1}{\varepsilon}D\Psi(\boldsymbol c)-D_{c}W(\boldsymbol c,{\mathcal E}(\b u))\\
        -\nabla\cdot D_{{\mathcal E}}W(\boldsymbol c,{\mathcal E}(\b u))\\         
    \end {array} 
\right).
\end{align}
We denote by $D_c$ and $D_{\mathcal E}$ the differentials with respect to $\b c$ and ${\mathcal E}$, respectively. In the case of a nonsmooth obstacle potential, $\Psi$ is given as the sum of a differentiable and a non-differentiable convex function and the derivative $D\Psi(\b c)$ has to be understood as the sum of the differentiable part plus the subdifferential of the non-differentiable convex summand, and so the first component of $\eqref{elastic-Allen}$ will result in a variational inequality, see Section \ref{Obstacle}. We have $D_c W(\b c,{\mathcal E})=-{\mathcal E}^*:{\mathcal C}({\mathcal E}-{\mathcal E}^*(\b c))$ and $D_{\mathcal E} W(\b c,{\mathcal E})={\mathcal C}({\mathcal E}-{\mathcal E}^*(\b c))$.\\

\noindent We assume now that a volume force $\b f$ acts on $\Omega_T:=\Omega\times(0,T)$ and a surface load $\b g\in L^2(0,T;L^2(\Gamma_g,\R^d))$ acts on $\Gamma_g\subset\Gamma:=\partial\Omega$ until a given time $T>0$. Then with $\Gamma_D:=\Gamma\setminus\Gamma_g$, $\Gamma_T:=\Gamma\times(0,T)$ and the outer unit normal $\b n$ the mechanical system is given by\vspace*{-0.3cm}

\begin{align}\label{state}
\left\{ \begin{array}{rcll}
-\nabla\cdot D_{{\mathcal E}}W(\boldsymbol c,{\mathcal E}(\b u)) & = & \boldsymbol 0 & \text{ in }\Omega,\\
{\boldsymbol u} & = & {\boldsymbol 0} & \text{ on }\Gamma_{D},\\
D_{{\mathcal E}}W(\boldsymbol c,{\mathcal E}(\b u))\cdot\boldsymbol n&= & {\boldsymbol g} & \text{ on }\Gamma_{g}
\end{array}\right.
\end{align}   
 which has to hold for a.e. $t\in (0,T)$, and the Allen-Cahn system is given by 
 \begin{align}\label{Allen-Cahn}
\left\{ \begin{array}{rcll}
\varepsilon\partial_{t}\boldsymbol c-\varepsilon\Delta\boldsymbol c+\frac{1}{\varepsilon}D\Psi(\boldsymbol c)+D_{\boldsymbol c}W(\boldsymbol c,{\mathcal E}(\b u)) & = & \boldsymbol f & \text{ in }\Omega_T,\\
\nabla{\boldsymbol c}\cdot\boldsymbol n & = & {\boldsymbol 0} & \text{ on }\Gamma_T,\\
{\boldsymbol c}(0) & = & {\boldsymbol c}_0 & \text{ in }\Omega
\end{array}\right.
\end{align} 
\noindent in case of a smooth potential $\Psi$. Our aim in this paper is to transform an initial phase distribution $\boldsymbol c_{0}:\Omega\rightarrow\mathbb R^N$ with minimal cost of the controls to some desired phase pattern $\boldsymbol c_{T}\in\b L^2(\Omega):=L^2(\Omega,\mathbb R^N)$ at a given final time $T>0$ while tracking a desired evolution $\boldsymbol c_d\in\boldsymbol L^2(\Omega_T):=L^2(0,T;\b L^2(\Omega))$. Hence we consider the following objective functional:

\vspace*{-0.3cm}
\begin{align}\label{tracking-type}
J(\boldsymbol c,\boldsymbol f,\boldsymbol g):=&\frac{\nu_{T}}{2}\|\boldsymbol c(T,\cdot)-\boldsymbol c_{T}\|^{2}_{\boldsymbol L^{2}(\Omega)}+\frac{\nu_{d}}{2}\|\boldsymbol c-\boldsymbol c_{d}\|^{2}_{\b L^{2}(\Omega_T)}+\nonumber\\
&+\frac{\nu_{f}}{2\varepsilon}\|\boldsymbol f\|^{2}_{\b L^2(\Omega_T)}+
\frac{\nu_{g}}{2}\|\boldsymbol g\|^{2}_{L^2(0,T;L^2(\Gamma_{g},\mathbb R^d))}.
\end{align}

\noindent This leads to the following optimal control problem:

\vspace*{-0.2cm}
\begin{align}\label{e:OverallUpperLevelProblemSmooth}
({\mathcal P})\:\left\lbrace\begin{array}{ll}
\min & J(\boldsymbol c,\boldsymbol f,\boldsymbol g)\\
\text{over }&(\boldsymbol c,\boldsymbol f,\boldsymbol g)\in \boldsymbol{\mathcal V}\times\b L^2(\Omega_T)\times L^2(0,T;L^2(\Gamma_{g},\mathbb R^d))\\
{\text{s.t.}} &(\ref{state})\text{ and }(\ref{Allen-Cahn})\text{ hold}
\end{array}\right.
\end{align}
with $\boldsymbol{\mathcal V}:=L^{\infty}(0,T;\boldsymbol H^1(\Omega))\cap H^1(0,T;\boldsymbol L^2(\Omega))\cap L^2(0,T;\boldsymbol H^2(\Omega))$. We assume, that the Dirichlet part $\Gamma_D$ has positive $(d-1)$-dimensional Hausdorff measure and introduce the notation $H^1_D(\Omega,\R^d):=\{\boldsymbol u\in H^1(\Omega,\mathbb R^d)\mid \boldsymbol u|_{\Gamma_{D}}=\b 0\}$. Later on we will use also the space $\boldsymbol{\mathcal W}(0,T):= L^2(0,T;\boldsymbol H^1(\Omega))\cap H^1(0,T;\boldsymbol H^1(\Omega)^*)$.

\section {Existence theory and first-order optimality conditions}
In this section we discuss the existence of a minimum and the derivation of first-order necessary optimality systems. First we present the smooth potential case. Here, the standard optimization theory in function spaces is applicable and delivers a first-order necessary optimality system. Afterwards, we focus on the control problem with an obstacle potential leading to an optimal control problem with variational inequalities. Hence this belongs to the class of MPECs, where the standard control theory is in general not applicable. Here we employ a penalty approach for the problem without distributed control and a relaxation approach for the model without elasticity.

\subsection{Smooth $\Psi$}\label{smooth}
We start by considering the setting without volume force, i.e. $\b f\equiv\b 0$. In a system with two phases, i.e. $N=2$, the problem can be reduced to a single unkown by defining $c:=c_1-c_2$, which results in a scalar problem. One typical choice of a smooth potential is then the double-well potential $\Psi(c)=\frac14(c^2-1)^2$. The scalar case with this $\Psi$ is studied extensively in \cite{Hecht} without tracking $c_d$, i.e. $\nu_d=0$. For a regularized obstacle potential $\Psi_\sigma$ (see Subsection~\ref{Pen}) the vector-valued case with possibly $\nu_d\neq0$ is discussed in \cite{FH}. However, $\Psi_\sigma$ is not a physical potential. The following theorem summarizes the results of \cite{FH,Hecht}. 
\begin{thm} Let $({\mathcal P})$ be given as a scalar problem for $N=2$ with potential $\Psi=\frac14\left(c^2-1\right)^2$ and $\nu_d=0$ or for $N\geq2$ and $\nu_d\geq0$ arbitrary with a regularized obstacle potential $\Psi_\sigma$ as mentioned above. For fixed initial distribution $\boldsymbol c_0\in \boldsymbol H^1(\Omega)$ and given surface load $\boldsymbol g\in L^2(0,T;L^2(\Gamma_g,\mathbb R^d))$ there exists a unique solution $(\boldsymbol c,\boldsymbol u)\in \boldsymbol{\mathcal V}\times L^2(0,T;H^1_D(\Omega,\R^d))$ of $\eqref{state}$-$\eqref{Allen-Cahn}$ and hence the solution operator $\boldsymbol S:L^2(0,T;L^2(\Gamma_g,\R^d))\rightarrow \boldsymbol{\mathcal V}\times L^2(0,T;H^1_D(\Omega,\R^d))$ with its components $\boldsymbol S(\boldsymbol g):=(\boldsymbol S_{1}(\boldsymbol g),\boldsymbol S_{2}(\boldsymbol g))=(\b c,\b u)$ is well-defined. 
\end{thm}
\noindent Then the control problem $({\mathcal P})$ is equivalent to minimizing the reduced cost functional $j(\boldsymbol g):=J(\boldsymbol S_{1}(\boldsymbol g),\boldsymbol g)$ over $L^2(0,T;L^2(\Gamma_g,\R^d))$. This result is established by applying energy methods to a time-discretized version of $\eqref{state}$-$\eqref{Allen-Cahn}$  and showing a series of uniform a priori estimates for the time discretized solutions, where one has to consider the particular functions $\Psi$ and $\Psi_\sigma$, respectively, and the coupling of the systems. By the direct method in the calculus of variations one can then show existence of a minimizer for $({\mathcal P})$. The differentiability of the solution operator can be shown by an implicit function argument and thus we can differentiate the reduced cost functional to obtain the following necessary optimality condition:
\begin{thm}
	Every minimizer $\b g\in L^2(0,T;L^2(\Gamma_g,\R^d))$ of $j$ fulfills the following optimality system:  $\eqref{state}$, $\eqref{Allen-Cahn}$ and\vspace*{-0.15cm}
\begin{align}\label{GradientEquationBoundary}\boldsymbol q+\nu_g\boldsymbol g=\boldsymbol 0\quad\text{ a.e. on }(0,T)\times\Gamma_g,\end{align} 
\begin{align}\label{e:NecessaryOptimalitySystem_Eqp1}
\left\{ \begin{array}{rcll}
-\varepsilon\partial_t\b p-\varepsilon\Delta\b p+\frac{1}{\varepsilon}D^2\Psi(\boldsymbol c)\b p+D_{p}W(\b p,\Epsilon(\b q))) & = & \nu_d(\boldsymbol c-\boldsymbol c_d) & \text{ in }\Omega_T,\\
\nabla{\boldsymbol p}\cdot\boldsymbol n & = & {\boldsymbol 0} & \text{ on }\Gamma_T,\\
\varepsilon{\boldsymbol p}(T) & = & \nu_T(\boldsymbol c(T)-\boldsymbol c_T) & \text{ in }\Omega,
\end{array}\right.
\end{align} 
\vspace*{-0.2cm}
\begin{align}\label{e:NecessaryOptimalitySystem_Eqp2}
\left\{\begin{array}{rcll}
-\nabla\cdot D_{{\mathcal E}}W(\boldsymbol p,{\mathcal E}(\b q)) & = & \boldsymbol 0 & \text{ in }\Omega,\\
{\boldsymbol q} & = & {\boldsymbol 0} & \text{ on }\Gamma_{D},\\
D_{{\mathcal E}}W(\boldsymbol p,{\mathcal E}(\b q))\cdot\boldsymbol n&= & {\boldsymbol 0} & \text{ on }\Gamma_{g}.
\end{array}\right.
\end{align} 
\end{thm}

\noindent For a setting without elasticity but with distributed control, i.e. $\boldsymbol f\not\equiv\boldsymbol 0$ and arbitrary $\nu_d,\nu_T\geq 0$, we refer for instance to \cite{FS1}. There, the scalar case, i.e. $N=2$ as above, is considered with a penalized double obstacle potential $\Psi_\sigma$. Moreover, the optimality system is investigated rigorously and is given by $\eqref{Allen-Cahn}$, $\eqref{e:NecessaryOptimalitySystem_Eqp1}$ without elastic energy together with the gradient equation
\begin{align}\label{GradientEquation}\b p+\frac{\nu_f}{\varepsilon}\b f=\b 0\quad\text{ a.e. in }\Omega_T.\end{align}

\subsection{Obstacle potential}\label{Obstacle}

In the case of an obstacle potential each component of $\boldsymbol c$ stands, in contrast to the smooth potential, exactly for the fraction of one phase. Hence the phase space is the Gibbs simplex $\boldsymbol G:=\{\boldsymbol v\in\mathbb R^N\mid v_i\geq 0,\sum_{i=1}^N v_i=1\}$ and the bulk potential $\Psi:\mathbb R^N\rightarrow\mathbb R\cup\{\infty\}$ is the multi-obstacle potential $\Psi(\boldsymbol v):=\Psi_{0}(\boldsymbol v)+I_{\boldsymbol G}(\boldsymbol v)$, where e.g. $\Psi_{0}(\boldsymbol v):=-\frac12\|\b v\|^2$, which we consider, and $I_{\boldsymbol G}$ is the indicator function of the Gibbs simplex. The differential of the indicator function has to be understood in the sense of subdifferentials, and thus the Allen-Cahn system $\eqref{Allen-Cahn}$ results in a variational inequality, which can also be written in the following form (see \cite{Garcke-Allen-Cahn-2}):
\vspace*{-0.05cm}
\begin{align}\label{Allen-Cahn-VI}
\left\{ \begin{array}{rcll}
\varepsilon\partial_{t}\boldsymbol c-\varepsilon\Delta\boldsymbol c-\boldsymbol P_{\Sigma}\left(\frac{1}{\varepsilon}(\boldsymbol c+\boldsymbol\xi)-D_{\boldsymbol c}W(\boldsymbol c,{\mathcal E}(\b u))\right) & = & \boldsymbol f & \text{ in }\Omega_T,\\
\nabla{\boldsymbol c}\cdot\boldsymbol n & = & {\boldsymbol 0} & \text{ on }\Gamma_T,\\
{\boldsymbol c}(0) & = & {\boldsymbol c}_0 & \text{ in }\Omega,
\end{array}\right.
\end{align} 
together with the complementarity conditions
\begin{align}\label{CC-all}
&\boldsymbol c\geq\boldsymbol 0 \text{ a.e. in } \Omega_{T},\: \boldsymbol \xi\geq\boldsymbol 0 \text{ a.e. in } \Omega_{T},\: (\boldsymbol \xi,\boldsymbol c)_{\boldsymbol L^2(\Omega_T)}=0,
\end{align}  
the additional constraint $\b c\in\b\Sigma:=\{\b v\in\R^N\mid\sum_{i=1}^Nv_i=1\}$ a.e. in $\Omega_T$ and the requirement $\b f\in\b{T\Sigma}:=\{\b v\in\R^N\mid\sum_{i=1}^Nv_i=0\}$ a.e. in $\Omega_T$. Here $\boldsymbol P_{\Sigma}:\mathbb R^N\rightarrow\boldsymbol{T\Sigma}$ is the projection operator defined by $\boldsymbol P_{\Sigma}\boldsymbol v:=\boldsymbol v-\boldsymbol 1\frac{1}{N}\sum\limits_{i=1}^{N}v_{i}$. The variable $\b\xi$ can be interpreted as a Lagrange multiplier corresponding to the constraint $\boldsymbol c\geq \boldsymbol 0$, and as a slack variable used for reformulating the variational inequality into a standard MPEC problem. Denoting $\b L^2_{\b{T\Sigma}}\left(\Omega_T\right):=\left\{\b v\in\b L^2(\Omega_T)\mid\b v\in\b{T\Sigma}\text{ a.e. in }\Omega_T\right\}$ and $\boldsymbol {\mathcal V}_{\b{T\Sigma}}$, $\boldsymbol {\mathcal V}_{\b{\Sigma}}$ respectively, the optimal control problem in the case of the obstacle potential is given by
\begin{align}\label{e:OverallUpperLevelProblemObstacle}
({\mathcal P}_{0})\:\left\lbrace\begin{array}{ll}
\min & J(\boldsymbol c,\boldsymbol f,\boldsymbol g)\\
\text{over }&(\boldsymbol c,\boldsymbol f,\boldsymbol g)\in \boldsymbol{\mathcal V}_{\boldsymbol\Sigma}\times\b L^2_{\boldsymbol{T\Sigma}}(\Omega_T)\times L^2(0,T;L^2(\Gamma_{g},\mathbb R^d))\\
{\text{s.t.}} &(\ref{state}), (\ref{Allen-Cahn-VI})\text{ and }(\ref{CC-all})\text{ hold}.
\end{array}\right.
\end{align}
The optimization problem $({\mathcal{P}}_{0})$ belongs to the problem class of so-called MPECs (Mathematical Programs with Equilibrium Constraints) which violate classical NLP constraint qualifications. In the next two subsections we present results concerning first-order necessary optimality systems obtained by the penalization approach, see \cite{FH}, or the relaxation approach, see \cite{FS2}. These techniques have been discussed also in \cite{Berg,HK, HW}.
\subsubsection{Penalization approach without distributed control}\label{Pen}
In this section we discuss the penalization approach for the case $\boldsymbol f\equiv\boldsymbol 0$. For the scalar Allen-Cahn case with $\boldsymbol f\not\equiv\boldsymbol 0$ but without elasticity we refer the reader to \cite{FS1}. Following \cite{FH} we replace the indicator function for the Gibbs simplex by a convex function $\tilde\psi_{\sigma}\in C^2(\mathbb R)$, $\sigma\in(0,\frac14)$, given by $\tilde\psi_{\sigma}(r):=0$ for $r\geq 0$, $\tilde\psi_{\sigma}(r):=-\frac{1}{6\sigma^2}r^{3}$ for $-\sigma<r<0$ and $\tilde\psi_{\sigma}(r):=\frac{1}{2\sigma}\left(r+\frac{\sigma}{2}\right)^{2}+\frac{\sigma}{24}$ for $r\leq -\sigma$, and define the regularized potential function by $\Psi_{\sigma}(\boldsymbol c)=\Psi_{0}(\boldsymbol c)+\hat{\Psi}(\boldsymbol c)$ with $\hat{\Psi}(\boldsymbol c):=\sum\limits_{i=1}^{N}\tilde\psi_{\sigma}(c_{i})$. For the resulting penalized optimal control problem denoted by $({\mathcal P}_{\sigma})$, exploiting techniques as in Section \ref{smooth}, we derive for $\sigma\in(0,\frac14)$ first-order necessary optimality conditions. Proving a priori estimates, uniformly in $\sigma\in(0,\frac14)$,
employing compactness and monotonicity arguments and using the definition $\b{\mathcal 
W}_0(0,T)=\{\b v\in\b{\mathcal W}(0,T):\b v(0,\cdot)=\b 0\}$ with dual space $\b{\mathcal 
W}_0(0,T)^*$, we are able to show the following existence and approximation result:
\begin{thm}\label{PenApproachExist} Whenever $\{\boldsymbol g_{\sigma}\}\subset L^2(0,T;L^2(\Gamma_g,\mathbb R^d))$ is a sequence of optimal controls for $({\mathcal P}_{\sigma})$ with the sequence of corresponding states $(\boldsymbol c_{\sigma},\boldsymbol u_{\sigma},\boldsymbol\xi_{\sigma})\in \boldsymbol{\mathcal V}_{\b{\Sigma}}\times L^{2}(0,T;H^1_D(\Omega,\R^d))\times \boldsymbol L^{2}(\Omega_{T})$, where $-\boldsymbol\xi_\sigma:=D\hat{\Psi}(\boldsymbol c_\sigma)$, and adjoint variables $(\b p_{\sigma},\b q_{\sigma},\b\zeta_{\sigma})\in \boldsymbol{\mathcal V}_{\boldsymbol{T\Sigma}}\times L^{2}(0,T;H^1_D(\Omega,\R^d))\times \boldsymbol L^{2}(\Omega_{T})$, where $-\boldsymbol\zeta_\sigma:=D^2\hat{\Psi}(\boldsymbol c_\sigma)\b p_\sigma$, there exists a subsequence, which is denoted again by $\{\b g_\sigma\}$, that converges weakly to $\b g$ in $L^2(0,T;L^2(\Gamma_g,\R^d))$. Moreover, $\b g$ is an optimal control of
$({\mathcal P}_{0})$ with corresponding states $(\boldsymbol c,\boldsymbol u,\boldsymbol\xi)\in \boldsymbol{\mathcal V}_{\b{\Sigma}}\times \boldsymbol L^{2}(\Omega_{T})\times L^{2}(0,T;H^1_D(\Omega,\R^d))$ and adjoint variables $(\b p,\b q,\b\zeta)\in L^2(0,T;\b H^1(\Omega))\times\linebreak[4] L^2(0,T;H^1_D(\Omega,\R^d))\times \b{\mathcal W}_0(0,T)^*$ and we have for $\sigma\searrow0$:
\begin{align}\label{convergences}
\left .\begin{array}{llllll}
\boldsymbol c_{\sigma}&\longrightarrow \boldsymbol c&\mbox{ weakly }&\text{in}&H^1(0,T;\b L^2(\Omega))\cap L^2(0,T;\b H^2(\Omega)),\\
\boldsymbol u_{\sigma}&\longrightarrow \boldsymbol u&\mbox{ weakly }&\text{in}&L^{2}(0,T;H^1_D(\Omega,\R^d)),\\
\boldsymbol\xi_{\sigma}&\longrightarrow\boldsymbol\xi&\mbox{ weakly }&\text{in}&\boldsymbol L^{2}(\Omega_{T}),\\
{\b p}_{\sigma}&\longrightarrow\b p&\mbox{ weakly 
}&\text{in}& L^2(0,T;\b H^1(\Omega)),\\
{\b q}_{\sigma}&\longrightarrow\b q&\mbox{ weakly }&\text{in}& 
L^2(0,T;H^1_D(\Omega,\R^d)),\\
\b P_{\b\Sigma}(\b\zeta_{\sigma})&\longrightarrow\b\zeta&\mbox{ weakly-star }&\text{in}& 
\b{\mathcal W}_0(0,T)^*.
\end{array} \raisetag{3\baselineskip} \right.
\end{align}
\end{thm}

\noindent Furthermore we obtain first order conditions:
\begin{thm}\label{PenalizationFirstOrderCond}
	The following optimality system holds for the limit elements $(\b g,\b c,\b u,\b\xi)$ with adjoint variables $(\b p,\b q,\b\zeta)$ of Theorem~\ref{PenApproachExist}:\\ $\eqref{state}$, $\eqref{GradientEquationBoundary}$, $\eqref{e:NecessaryOptimalitySystem_Eqp2}$, $\eqref{Allen-Cahn-VI}$, $\eqref{CC-all}$, $\b c\in\b\Sigma$, $\b f\in\b{T\Sigma}$ a.e. in $\Omega_T$ and  
\begin{align}
&-\frac{1}{\varepsilon}\b\zeta(\b v)+\varepsilon\int_0^T\left\langle\partial_t\b v,\b 
p\right\rangle\dt +\varepsilon\int_0^T\int_\Omega\nabla\b p\cdot\nabla\b v\dx\dt+\nonumber\\
&-\frac{1}{\varepsilon}\int_0^T\int_\Omega\b p\cdot\b v\dx\dt+\int_0^T\int_\Omega\b P_{\b\Sigma}(D_p W(\b p,\Epsilon(\b 
q)))\cdot\b v\dx\dt+\nonumber\\
&-\int_0^T\int_\Omega\nu_d({\boldsymbol c}-\boldsymbol c_d)\cdot\b v\dx\dt-\int_\Omega\nu_T({\b c}(T,\cdot)-\b c_T)\cdot\b v(T)\dx=0,\label{e:NecessaryOptimalitySystemLimitcp}
\end{align}
which has to hold for all $\b v\in\b{\mathcal W}_0(0,T)$. Moreover, the limit elements satisfy some sort of complementarity slackness 
conditions:
\begin{align}
&\lim_{\sigma\searrow 0}(\b\zeta_{\sigma},{\b 
p}_{\sigma})_{\b L^2\left(\Omega_T\right)}\leq 0,\label{e:ComplemenLimitSystem_zeta_p1}\\
&\lim_{\sigma\searrow 0}(\b\zeta_{\sigma},\max(\b 0,{\b c}_{\sigma}))_{\b L^2(\Omega_T)}=0,\label{e:ComplemenLimitSystem_zeta_c}\\
&\lim_{\sigma\searrow 0}({\b p}_{\sigma},{\b 
\xi}_{\sigma})_{\b L^2(\Omega_T)}=0.\label{e:ComplemenLimitSystem_p1}
\end{align}
\end{thm}

\subsubsection{Relaxation approach with distributed control and without elasticity}\label{relaxHassan}
Studying the control problem with distributed control, i.e. $\b f\not\equiv\b0$ in general, and without elasticity we use a relaxation approach. Details for our presented results can be found in \cite{FS2}. After reformulating as in $\eqref{Allen-Cahn-VI}-\eqref{CC-all}$ the Allen-Cahn system with the help of a slack variable $\b\xi$ into an MPEC, we add to the problem $({\mathcal P}_{0})$ an additional constraint $\frac{1}{2}\|\boldsymbol\xi\|^2_{\b L^2(\Omega_T)}\linebreak[3]\leq R$ and denote this modified optimization problem by $({\mathcal P}_{R})$. The constant $R$ is sufficiently large. This approach is also used in \cite{Berg} where the control of an obstacle problem is considered. As a first step we treat the state constraint $\b c\geq\b 0$, which usually raises problems concerning regularity, by adding a regularization term to $J$. I.e. we define $J_\gamma(\b c,\b f)=J(\b c,\b f)+\frac{1}{2\gamma\varepsilon}\sum\limits_{i=1}^N\|\max(0,\overline{\lambda}-\gamma c_{i})\|^{2}_{L^{2}( \Omega_{T})}$ where $\overline\lambda\in L^2(\Omega_T)$ is fixed, nonnegative and corresponds to a regular version of the multiplier associated to $\b c\geq\b 0$. Next we relax the complementarity condition to $(\b\xi,\b c)_{\b L^2(\Omega_T)}\leq\varepsilon\alpha_\gamma$ for some $\alpha_\gamma>0$. We denote this regularized relaxed version of $({\mathcal P}_{R})$ as $({\mathcal P}_{R,\gamma})$. Subsequently we are interested in $\gamma\nearrow\infty$ where simultaneously $\alpha_\gamma\searrow0$. We are able to use techniques from mathematical programming in Banach spaces, see \cite{Zowe}, and get an optimality system for $({\mathcal P}_{R,\gamma})$, where $\gamma$ is fixed. Considering $\gamma\nearrow\infty$ we then obtain optimality conditions for problem $({\mathcal P}_{R})$. Similar to the process in Section \ref{Pen} we have: for any $\gamma>0$ there exists a minimizer $(\b c_\gamma,\b f_\gamma,\b\xi_\gamma)\in\b V_{\b\Sigma}\times\b L^2(\Omega_T)\times\b L^2(\Omega_T)$ of $({\mathcal P}_{R,\gamma})$ with corresponding adjoint variables. Using the Lagrange multiplier $r_\gamma\in\R$ of the constraint $(\b\xi_\gamma,\b c_\gamma)_{\b L^2(\Omega_T)}\leq\varepsilon\alpha_\gamma$ one defines $\zeta_{\gamma,i}:=r_\gamma\xi_{\gamma,i}-\max(0,\overline\lambda-\gamma c_{\gamma,i})$ and $\b\zeta_\gamma:=(\zeta_{\gamma,i})_{i=1}^N$. Then we obtain:
 \begin{thm}
 Whenever $\{\b f_\gamma\}$ is a sequence of optimal controls $({\mathcal P}_{R,\gamma})$ with the sequence of corresponding states $(\b c_\gamma,\b\xi_\gamma)$ and adjoint variables $(\b p_\gamma,\b\zeta_\gamma)$, there exists a subsequence, which is denoted the same, with $\b f_\gamma\to\b f$ weakly in $\b L^2(\Omega_T)$ and $\b\zeta_\gamma\to\b\zeta$ weakly-star in $\b{\mathcal W}_0(0,T)^\ast$ as $\gamma\nearrow\infty$. The convergence of the variables $\b c_\gamma$, $\b\xi_\gamma$ and $\b p_\gamma$ is as in $\eqref{convergences}$. These limits fulfill the corresponding optimality system for $({\mathcal P}_{R})$ as in Theorem~\ref{PenalizationFirstOrderCond} without elasticity system but with distributed control, i.e. $\eqref{GradientEquation}$, $\eqref{Allen-Cahn-VI}$, $\eqref{CC-all}$, $\eqref{e:NecessaryOptimalitySystemLimitcp}$, $\b c\in\b\Sigma$, $\b f\in\b{T\Sigma}$ a.e. in $\Omega_T$ and the limits with $(\b p_\gamma,\b \zeta_\gamma)$ satisfy the complementarity slackness conditions $\eqref{e:ComplemenLimitSystem_zeta_p1}$-$\eqref{e:ComplemenLimitSystem_p1}$ for $\gamma\nearrow\infty$ instead of $\sigma\searrow0$. In addition we have the constraint $\frac12\|\b\xi\|_{\b L^2(\Omega_T)}^2\leq R$.
 \end{thm}
 \noindent The last inequality is in practice inactive using $R$ large enough.
\vspace*{-0.6cm}
\section{Numerics}
\vspace*{-0.1cm}
In this section we neglect elastic effects, but study smooth as well as nonsmooth obstacle potentials with distributed control numerically. 
\vspace*{-0.2cm}\subsection{Smooth potential}\label{Christoph} 
\begin{figure}[htb]
\centering\hspace*{-2cm}
\begin{minipage}{.75\textwidth}
  \centering
	\includegraphics[width=3cm]{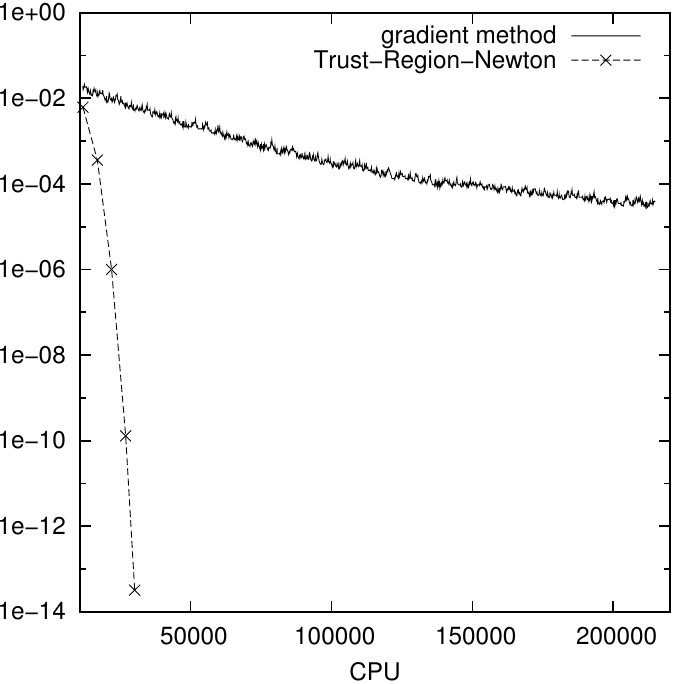}
	\caption{$\|\nabla j(f_k)\|_{L^2(\Omega_T)}$ depending on cpu-time for the TRN and the gradient method applied to have $c(T) = c_0$.}\label{fig:trgrad}
\end{minipage}\hspace*{-2cm}
\begin{minipage}{.6\textwidth}
  \centering
	\includegraphics[width=3cm]{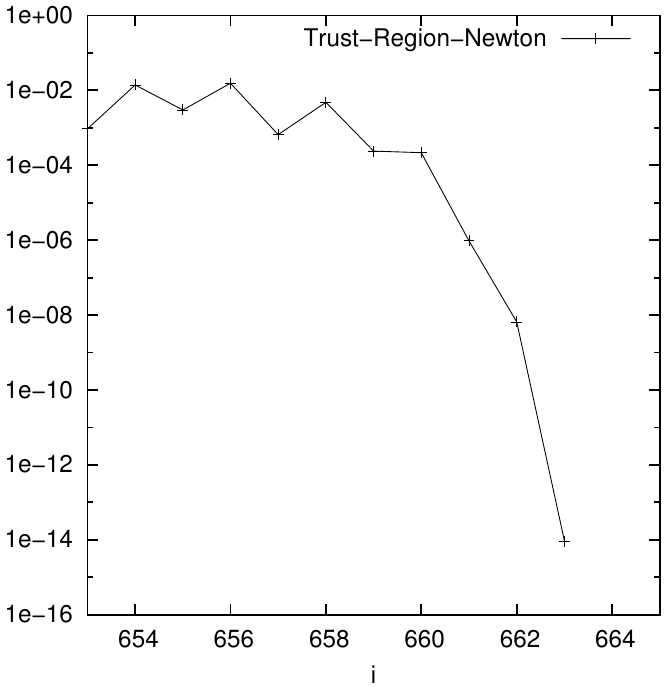}
	\caption{Newton residual for moving a vertical interface from left to right.}\label{fig:Residual-moving-line-short-time}
\end{minipage}\hspace*{-2cm}%
\end{figure}

{\bf Newton's method.} For smooth $\Psi$ we obtain an unconstrained optimal control problem when eliminating the state equation. Hence, numerical methods for unconstrained problems can be applied to the reduced problem\\
\centerline{$\min j(\u) := J(S(\u),\u),\quad \u\in \boldsymbol L^2(\Q)$.}
We choose the Trust-Region-Newton-Steihaug-cg (TRN) method, see \cite{Conn}, since it is capable of solving large scale optimization problems very efficiently because the underlying cg-solver is matrix-free and it attains the local convergence properties of Newton's method. Iteratively the model $m_k(\bdelta \u)= j(\u_k) + (\nabla j(\u_k),\bdelta \u)_{\boldsymbol L^2(\Omega_T)} + \frac{1}{2}(\nabla^2j(\u_k)\bdelta \u, \bdelta \u)_{\boldsymbol L^2(\Omega_T)}$ is minimized within a trust-region and the method is stopped if $\|\nabla j(\u_k)\|_{L^2(\Omega_T)} < tol$. \\
Based on Section \ref{smooth} the $L^2$-gradient is given by $\nabla j(\u)=\frac{\nu_f}{\varepsilon}\u + \p$. The Hessian we derived formally for $\nu_d=0$, see \cite{Rupprecht}, and is given by $\nabla^2j(\u)\bdelta \u = \frac{\nu_f}{\varepsilon}\bdelta \u + \bdelta \p$, where $\bdelta \p$ can be calculated by first solving the linear forward equation $\varepsilon\partial_t\bdelta \y -\varepsilon\Delta\bdelta \y + \frac{1}{\varepsilon}D^2\Psi(\y)\bdelta \y = \bdelta \u$  in $\Q$, $\nabla(\bdelta \y)\cdot\b n = \boldsymbol 0$ on $\Gamma_T$ and $\bdelta \y(0) = \boldsymbol 0$ in $\Omega$ and then solving the linear backward equation $-\varepsilon\partial_t\bdelta \p -\varepsilon\Delta\bdelta \p + \frac{1}{\varepsilon}D^2\Psi(\y)\bdelta \p = -\frac{1}{\varepsilon}D^3\Psi(\y)[\bdelta \y, \p,.]$ in $\Q$, $\nabla(\bdelta \p)\cdot\b n = \boldsymbol 0$ on $\Gamma_T$ and $\varepsilon\bdelta \p(T) = \nu_T\bdelta \y(T)$ in $\Omega$. The cost of one iteration of the algorithm consists in evaluating $j$, which means solving the nonlinear state equation, in calculating $\nabla j(\u)$, which means solving the linear adjoint equation, and in performing the Steihaug-cg method, where in each cg-iteration $\nabla^2j(\u)$ has to be evaluated in some direction $\bdelta\u$. For similar control problems gradient type methods have been used, see e.g. \cite{Hausser2010,Ohtsuka2012}. However, they cannot solve our problems in reasonable time.\\
The following numerical results summarize the investigations in \cite{Rupprecht}.\\
{\bf Discretization and error estimation.}
We consider an implicit and a semi-implicit Euler scheme in time. Although solving the semi-implicit discrete equations is much faster, it has the disadvantage that the two approaches ``first discretize then optimize'' and ``first optimize then discretize'' do not commute. This has been shown by looking upon the implicit discretization as a discontinuous Galerkin ansatz \cite{Rupprecht}. Thus we use semi-implicit discretization only in an initialization phase to compute an approximative optimal control, and use implicit discretization in the main phase.\\
In space we discretize with standard P1-elements. For equidistant meshes we implemented the TRN method with the toolbox FEniCS \cite{LoggMardalEtAl2012a}, exploiting the structure of the arising systems for equidistant meshes. The existing adaptive strategy for the Allen-Cahn equation without control uses a fine mesh on the interface and coarse mesh on the bulk regions, see e.g. \cite{Barrett}. However, with control, nucleation of a phase may appear. This cannot be resolved using the concept in \cite{Barrett}. Moreover, a method of adaptively controlling the time steps for Allen-Cahn equations is not available. Hence, for studying adaptive meshes we use the toolbox RoDoBo, where the TRN method together with a dual weighted residual (DWR) error estimator is implemented, see \cite{MeidnerVexler}. In our applications the DWR error estimator establishes both: adequate adaptive spatial meshes and adaptive time steps. For example in a nucleation situation the mesh in \cite{Barrett} is only fine when the new phase was already created, whereas the DWR mesh is also fine at timesteps before the nucleation process starts.\\
{\bf Numerical results.}
In all experiments we choose $d=N=2$, $\nu_T=1$, $\nu_d=0$, $\nu_f=0.01$, $\varepsilon=(14\pi)^{-1}$, $tol=10^{-13}$ and $\Omega=(-1,1)^2$. As mentioned above we reduce the problem to a scalar problem and use $\Psi(c) = \frac{1}{4}(c^2 - 1)^2$. Figure \ref{fig:trgrad} 
depicts the large speed up using the TRN method instead of the gradient method. Here the Newton residual $\|\nabla j(f_k)\|_{L^2(\Omega_T)}$ for the TRN method and the gradient method are listed for an example where $c(T)$ shall be the same circle as $c_0$. The cpu-time is still large for the TRN method using RoDoBo. However, using an equidistant mesh and therefore being able to exploit the structure of the problem, our implementation in FEniCS is significantly faster. Already the adjoint equation can be solved 25 times faster.\\
In order to get quadratic convergence of the Newton-cg method for smooth problems the inner tolerance $tol_{cg}$ has to be appropriate.
While one can decrease $tol_{cg}$ with the number of iterations, we set $tol_{cg}=10^{-13}$ in order to solve the inner problem nearly exact and the resulting numerical error does not influence the performance of the Newton-method.
In most experiments the Newton method converged just superlinearly which reveals that the problem is not smooth enough. Only in an experiment where a vertical interface is moved from left to right we could observe quadratic convergence in $\nabla j(f_k)$, see Figure \ref{fig:Residual-moving-line-short-time}. For the first 660 iterations in this example, an approximation of the model problem is computed with less than 40 Steihaug-cg-iterations. They always lie on the boundary of the trust-region. In the last three iterations the trust-region constraint stays inactive and then about 600 cg-iterations are necessary to solve the quadratic subproblem. In these last three outer iterations the convergence rate of Newton's method can be observed. Also in the other experiments in \cite{Rupprecht} the Steihaug-cg method performs only few inner cg-iterations when the trust-region constraint is active. In the last few steps the calculation of the unconstrained minimizer of $m_k$ is much more expensive.\\
Next we consider the situation where a circle in the center shall be split into two circles next to each other.
\begin{figure}[htb]
\centering
\subfloat{\includegraphics[width=2cm]{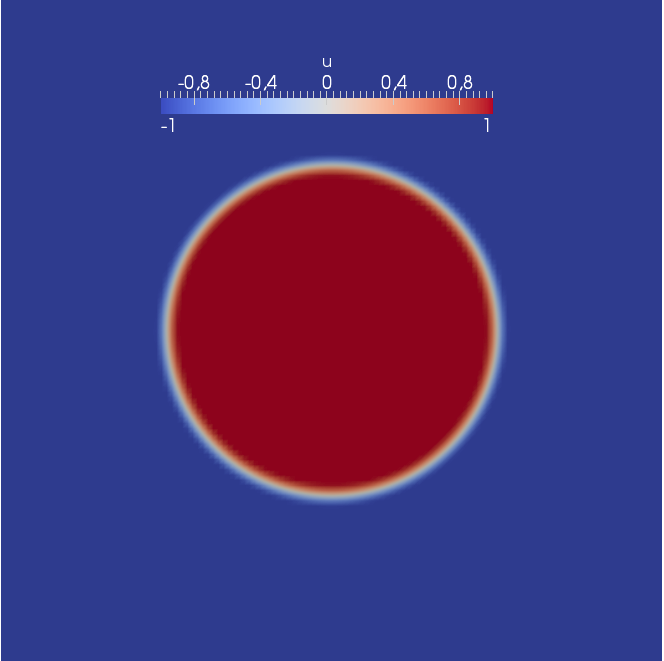}\label{fig:SCy0}}\hspace{.1cm}
\subfloat{\includegraphics[width=2cm]{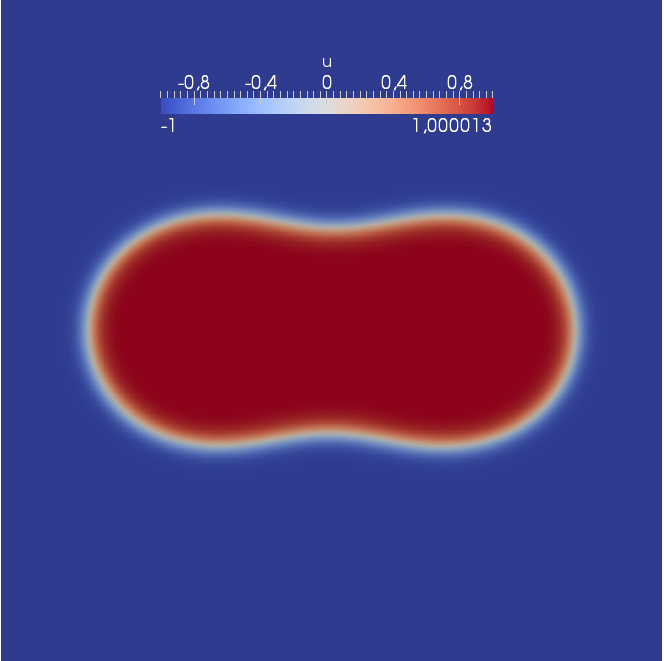}}\hspace{.1cm}
\subfloat{\includegraphics[width=2cm]{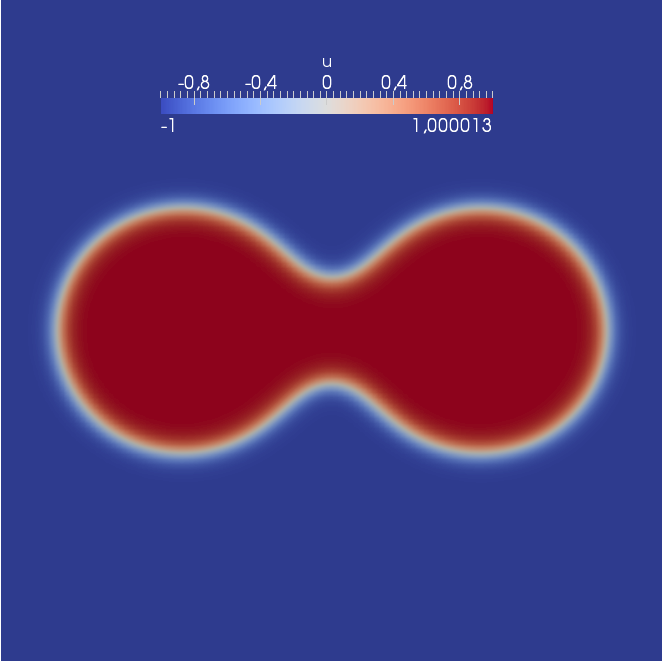}}\hspace{.1cm}
\subfloat{\includegraphics[width=2cm]{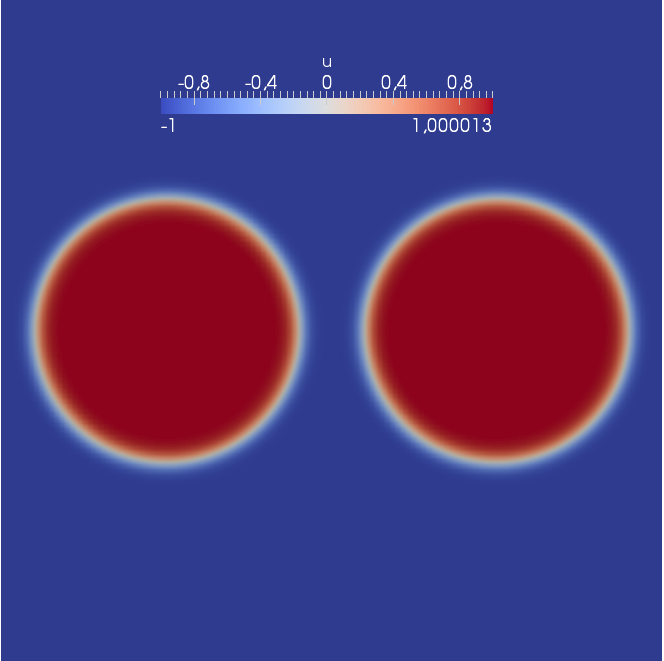}\label{fig:SCyT}}\hspace{.1cm}\\
\subfloat{\includegraphics[width=2cm]{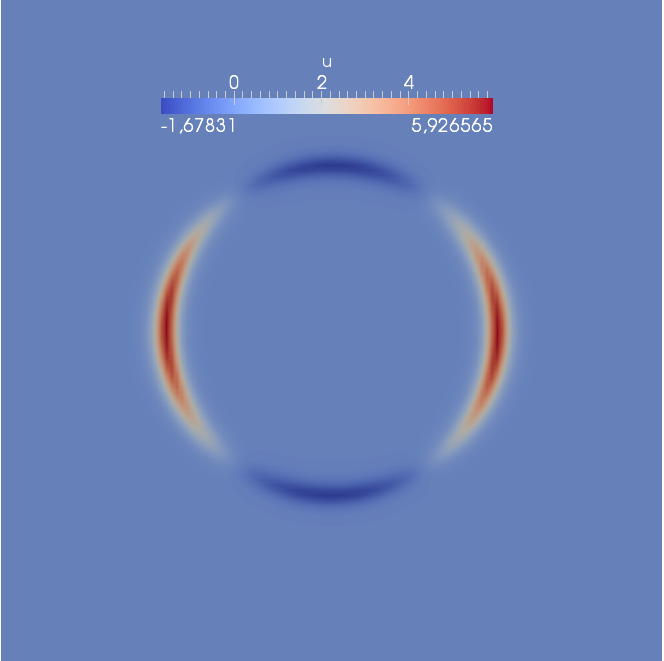}}\hspace{.1cm}
\subfloat{\includegraphics[width=2cm]{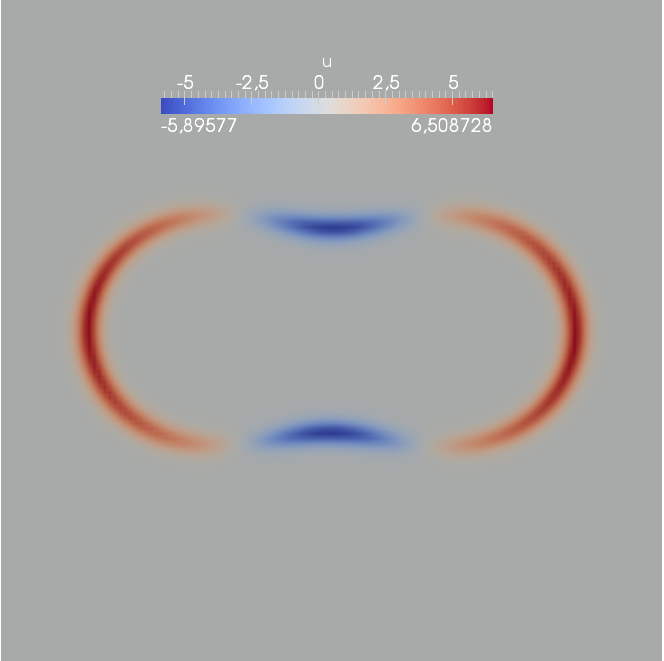}}\hspace{.1cm}
\subfloat{\includegraphics[width=2cm]{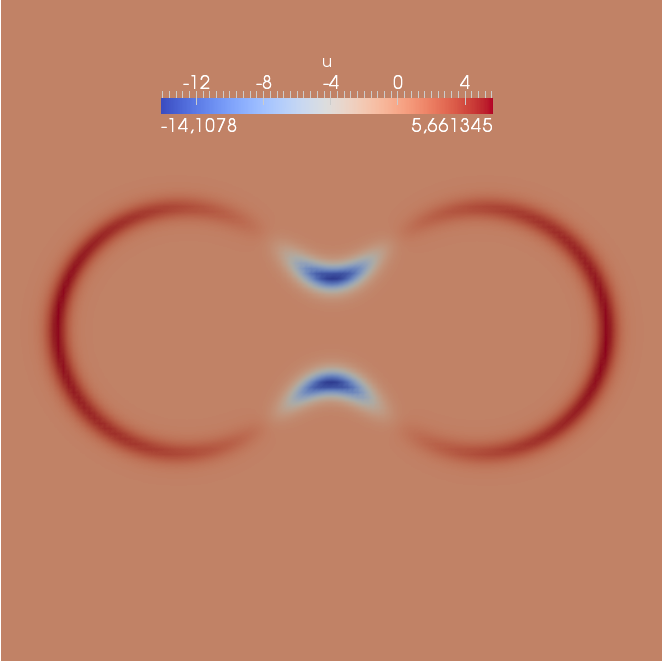}}\hspace{.1cm}
\subfloat{\includegraphics[width=2cm]{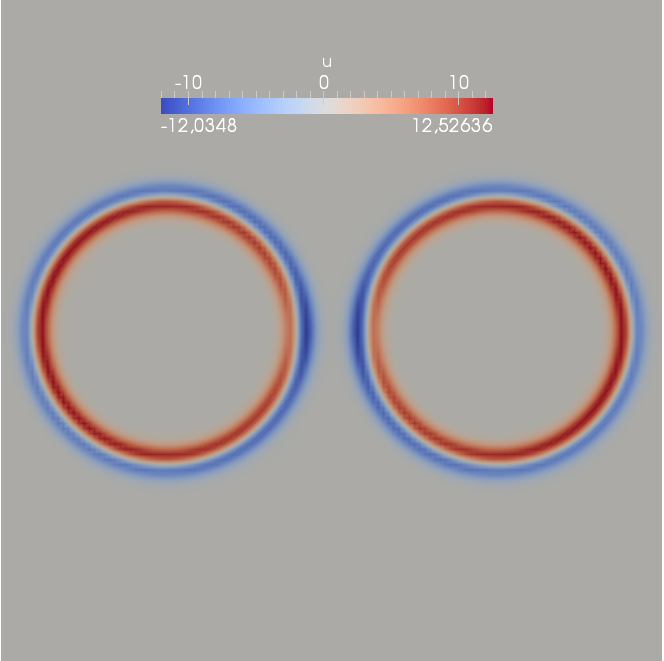}}
\caption{Optimal state (top) and optimal control (bottom) at times $t=0,\frac 1 2 T, \frac 3 4 T, T$, for a splitting circle scenario.}\label{fig:SCyu}
\vspace*{0.5cm}
\centering
\includegraphics[width=2.5cm]{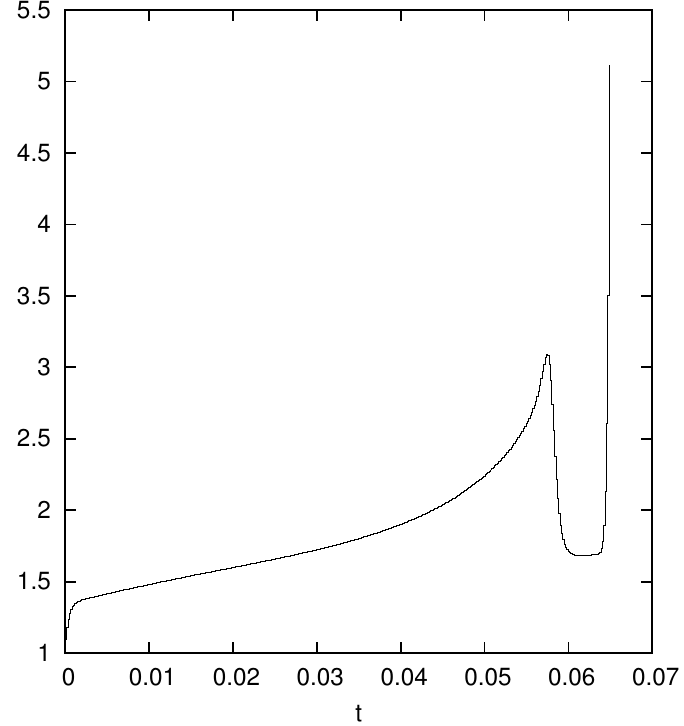}
\caption{$L^2(\Omega)$-norm of the control corresponding to Fig. \ref{fig:SCyu}.}\label{fig:SCL2u}
\end{figure}
Figure \ref{fig:SCyu} shows the optimal state and control. The circle is stretched horizontally until it separates into two circles. In Figure \ref{fig:SCL2u} the plot of $t\mapsto \|f(t)\|_{L^2(\Omega)}$ is depicted. The peak is at the time when the topological change occurs. The large increase of the cost at the end time is due to the fact that $c_T$ has a smaller interface thickness than proposed by the model with $\varepsilon$.\\
In the following we investigate the temporal mesh.
\begin{figure}[htb]
\centering
\includegraphics[width=8cm]{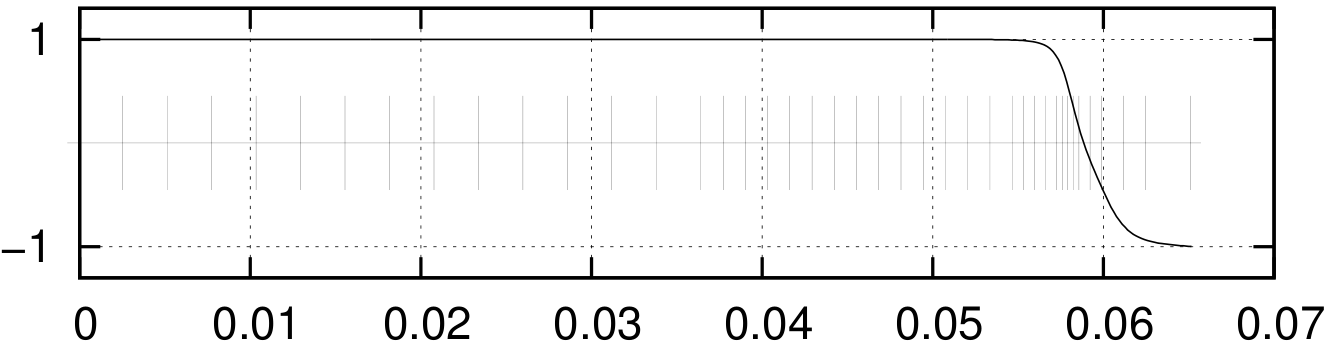}
\caption{Temporal mesh for a splitting circle scenario.}\label{fig:time-steps-topology-change}
\centering
\vspace*{0.5cm}
\includegraphics[width=8cm]{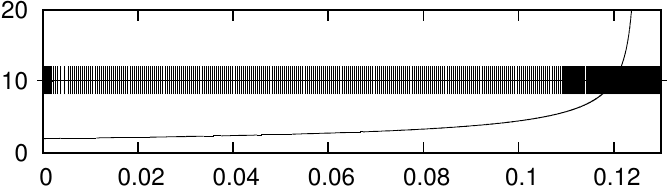}
\caption{Temporal mesh for the time evolution without control.}\label{fig:time-steps-interface-velocity}
\end{figure}
Figure \ref{fig:time-steps-topology-change} shows the time steps  created by the DWR error estimator together with the value $c(t)$ at the location $\b x=\boldsymbol 0$. We can see that the time steps are small before the pinching occurs, attain their minimum in the middle of the pinching process and are larger in the second half of the pinching process. To study how the time steps depend on the interface velocity we consider an experiment where a circle shrinks and vanishes at finite time. The end time is chosen in such a way that $f\equiv 0$ is the optimal control, i.e. the interface evolution is given by the Allen-Cahn equation without outer force. Figure \ref{fig:time-steps-interface-velocity} depicts the interface velocity together with the time steps. As expected, the larger the interface velocity becomes, the smaller the time steps have to be chosen.

\subsection{Obstacle potential}\label{Luise}
\begin{figure}[ht]
\hspace*{-7mm}\includegraphics[width=13cm]{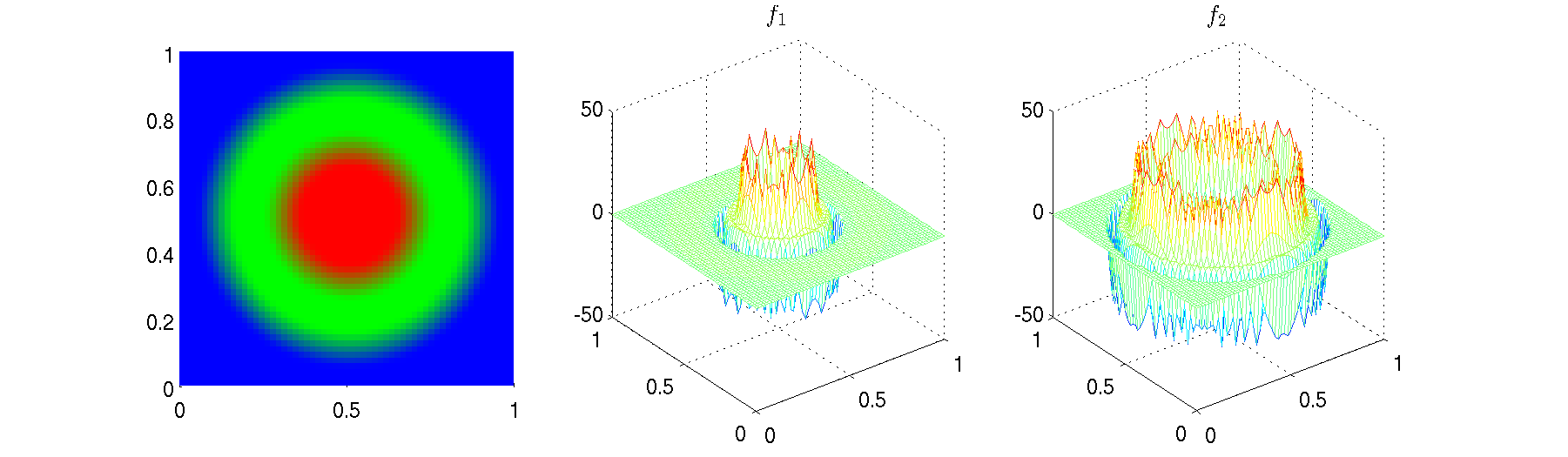}
\caption{\mbox{The state $\boldsymbol c$ and the control functions $f_1$ and $f_2$ for} three phases, which shall stay constant, at time $T=0.0003 $.}\label{contcirc}
\end{figure}

\begin{figure}[ht]\hspace*{-7mm}
\begin{tabular}{ccc}
\includegraphics[width=2.5cm]{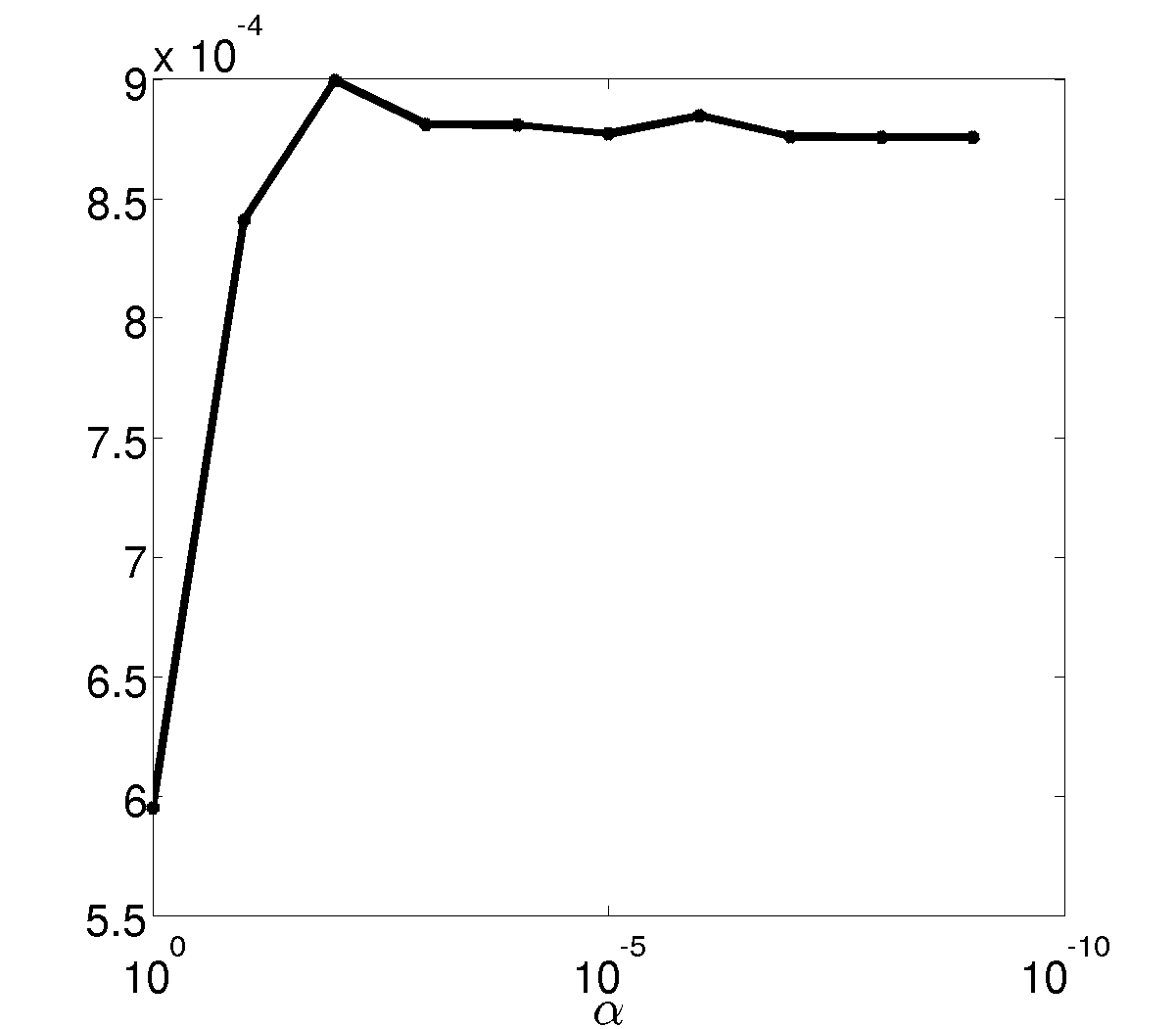}
&
\includegraphics[width=2.5cm]{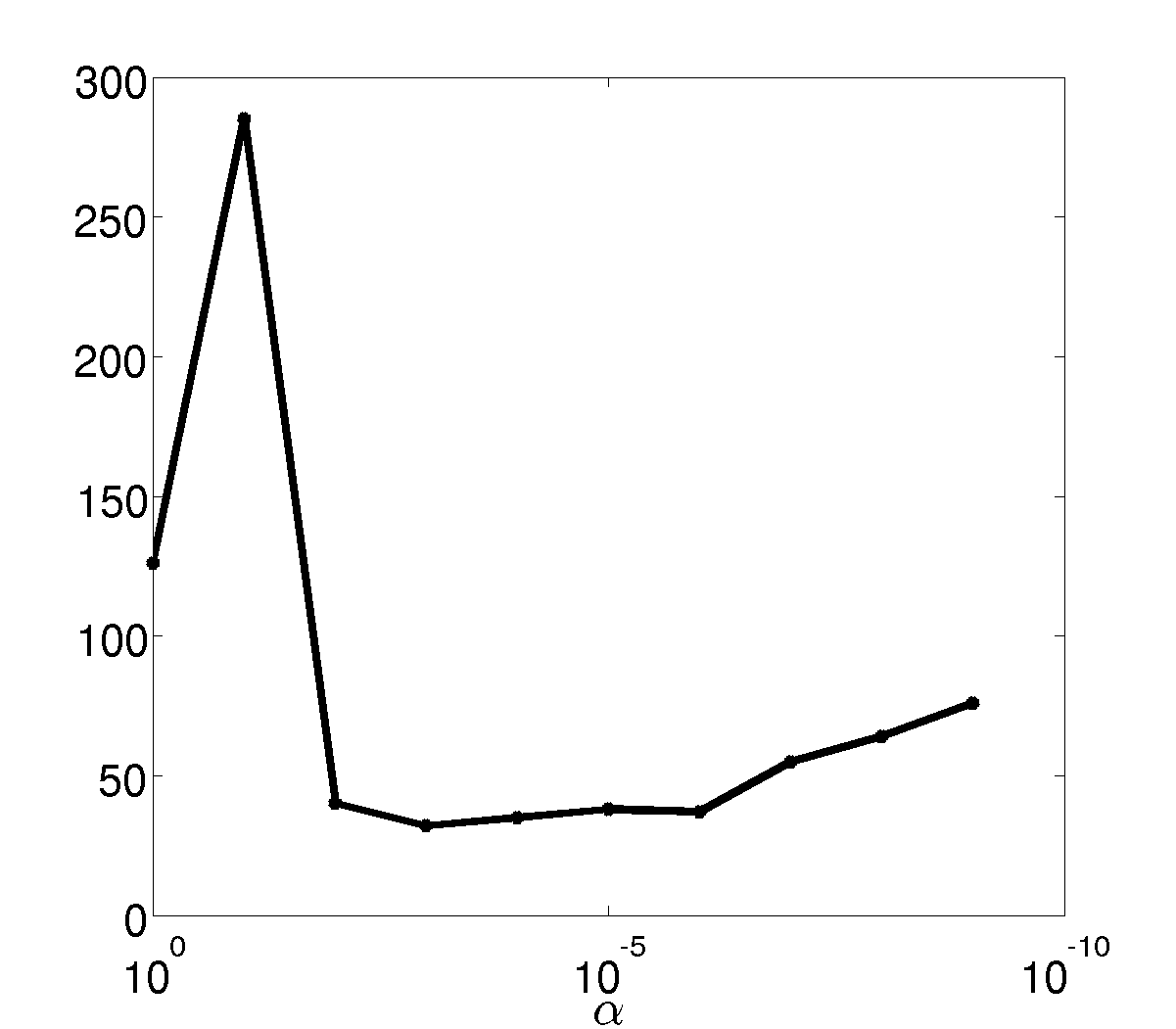}
&
\includegraphics[width=2.5cm]{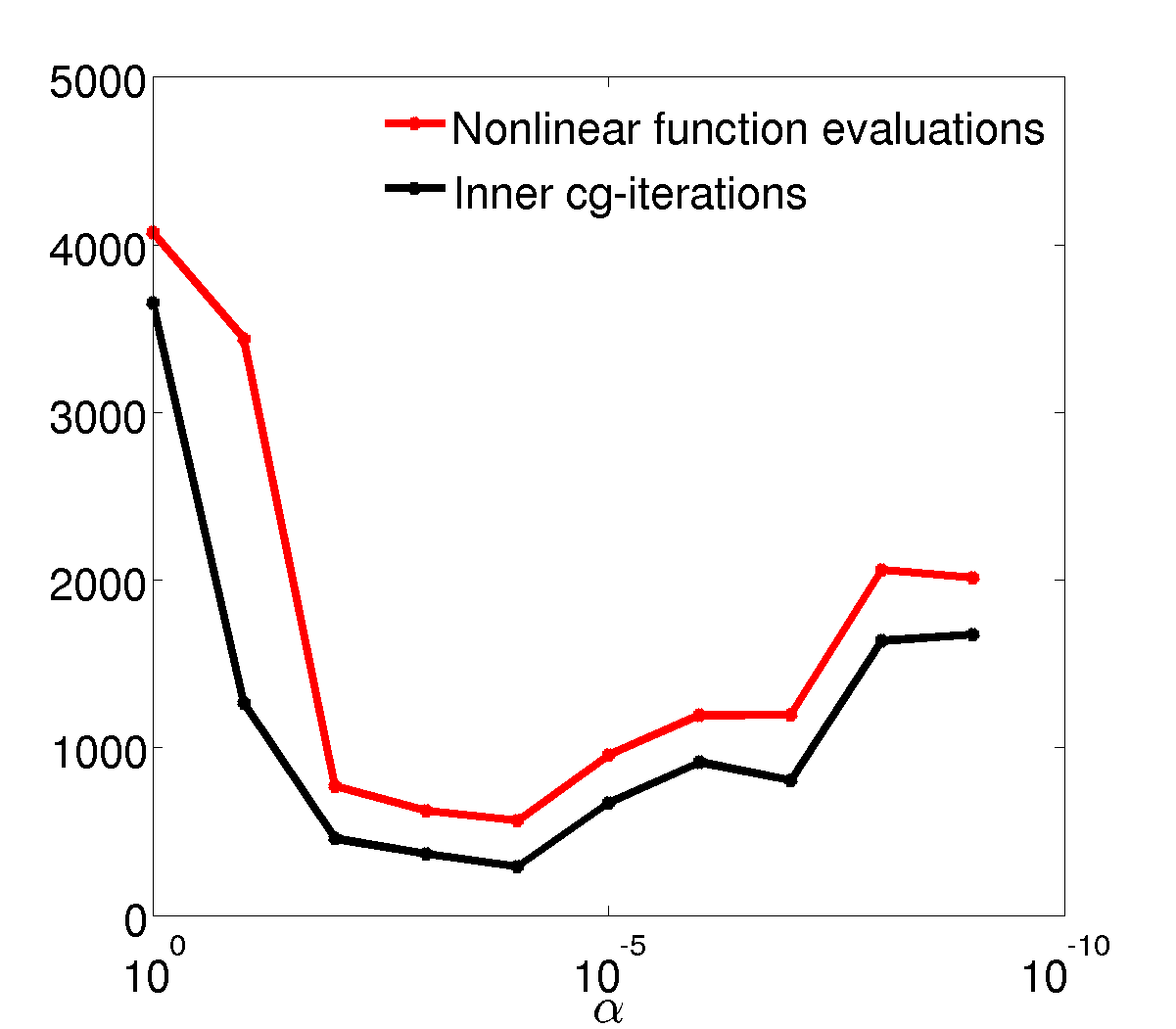}
\\
{\small Cost function $J$} & {\small Interior point} & {\small Nonlinear function evalua-} \\
&{\small iterations}&{\small tion and inner iterations}
 \end{tabular}
\caption{
Results for varying $\alpha$ for the example in Figure \ref{contcirc}.}\label{circiter}
\end{figure}

In the case of an obstacle potential we studied first the differences
in the approaches ``first discretize then optimize'' and
``first optimize then discretize''.
As in the smooth case the choice of discretization is essential.
We choose again an implicit discretization of the Allen-Cahn system,
which is understood as a discontinuous Galerkin discretization in time.
Hence the time integrals of functions are discretized by an iterated rectangle rule
using the right endpoints. This approximation is also used in the cost function.
We compared the discretized optimality system of the ansatz presented in Subsection \ref{relaxHassan} 
with the optimality system arising for  
the discretized optimization problem, where the first order conditions
(C-stationarity) are derived by the
relaxation approach in \cite{Scholtes} for finite dimensional MPECs 
assuming MPEC-LICQ.
In the latter only the complementarity condition is relaxed as in \ref{relaxHassan}
to  $(\b\xi_\alpha,\b c_\alpha)_{\b L^2(\Omega_T)}\leq \alpha$ .
The systems are identical apart from the additional constraint on 
$\boldsymbol\xi$, which is inactive in the numerics, and, as expected, the complementarity slackness conditions,
which hold pointwise for the ansatz, where the problem is discretized first
\cite{Michl}.\\
Our first numerical experiments are based on the MATLAB solver {\em fmincon} where
the discretized, relaxed optimization problem is solved  --- due to the memory limitations ---
using an interior point algorithm with internal cg--solver for
decreasing $\alpha$. 
The inital $\alpha_0=1$ is successively divided by ten and 
the solutions for $\alpha_i$ are used as initial data for the 
problem with relaxation parameter $\alpha_{i+1}$. 
In the first example with $N=3$ the goal is to keep the initial setting unchanged for the time interval $[0,0.0005]$,
where one phase in a circle is surrounded
by an annulus with a second phase and a third phase in the remainder of the domain
$\Omega =(0,1)^2$. Without any control the
two inner phases would vanish due to the curvature.
We set $\nu_T=1$, $\nu_d=10^4$, 
$\nu_f=0.001$, $\varepsilon = 0.1$ and the time step $\tau = 10^{-4}$ while the equidistant mesh size in space 
is $h= 1/59$. The phases stay nearly constant
as do the controls which we therefore list only for $T=0.0003$ and $\alpha=10^{-9}$ in Figure \ref{contcirc}.
The control $f_1$ is positive on the innermost interface to ensure that this circle does not shrink. However, noticeable is that $f_1$ is negative on the other 
interface, where it seems that $c_1$ would otherwise increase, i.e. phase one would develop.
In the same way $f_2$ is negative on the innermost circle while positive to hold the interface
constant on the outer circle. Correspondingly $f_3$ behaves. In Figure \ref{circiter} the first plot shows the values of the cost function $J$ neglecting the constant part for decreasing $\alpha$. 
For $\alpha \leq 10^{-3}$ it changes only mildly. 
The main effort of calculating the optimal control is used for large $\alpha$ as 
the other two plots in Figure \ref{circiter} indicate, which list the number of 
interior point iterations and the number of nonlinear function evaluations 
together with 
the cg-iterations.
They indicate also the expected cost if a more sophisticated implementation of
an optimization solver is employed.

\begin{figure}[ht]
\vspace*{-4mm}
\includegraphics[width=10cm]{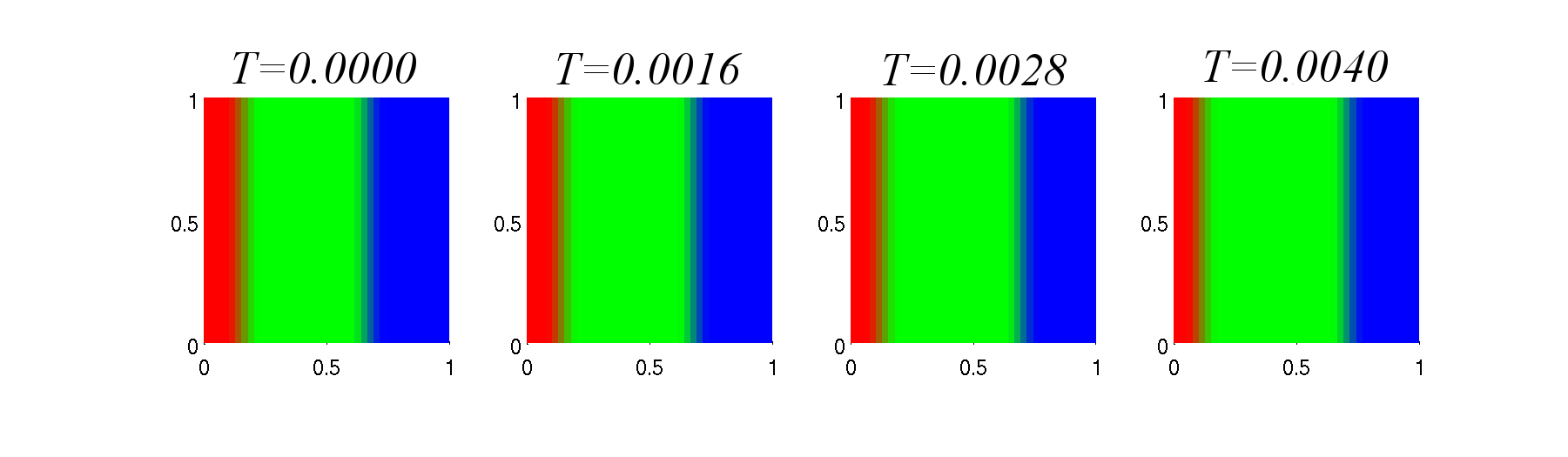}
\vspace*{-8mm}
\caption{\mbox{The state $\bf c$ for three phases 
for moving walls.}}\label{front}
\end{figure}

\begin{figure}[ht]
\hspace*{-17mm}\includegraphics[width=14cm]{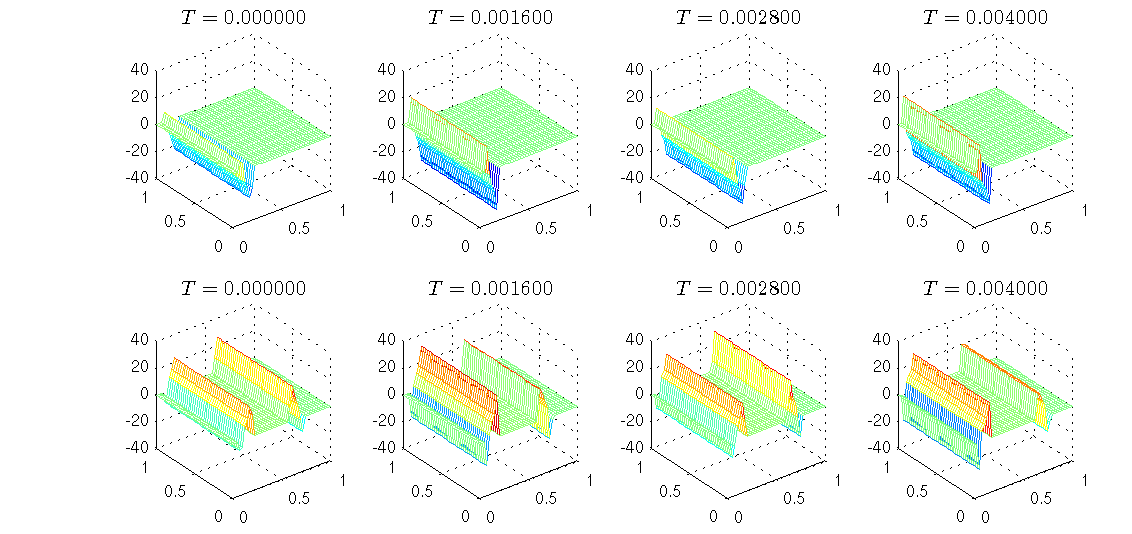}
\vspace*{-5mm}
\caption{Control functions $f_1$ in the first and $f_2$  in the second row corresponding to Figure \ref{front}.}\label{frontcont}
\end{figure}

\noindent In the next example three phases are vertically aligned. Since the
interfaces have no curvature the phases would stay constant without control.
However, in this experiment we set the target ${\bf c}_d$ such that in the 
end the enclosed phase  occupies a larger rectangle than the others
as the numerical result shows  in Figure \ref{front} for $\alpha = 10^{-9}$.
Hence the controls are now time dependent.
In the first row of Figure \ref{frontcont} $f_1$ is depicted and in the second $f_2$ while $f_3=-f_1-f_2$ is neglected.
As expected the controls work mainly on the interfaces. The control $f_2$ is positive at both
interfaces while the other two controls support the movement by negative force.
Like in the first example the value of the cost function
stays nearly constant for $\alpha \leq 10^{-3}$ and 
the substantial work of determining the optimal control is done for large $\alpha$. We therefore omit the figures.

\addcontentsline{toc}{chapter}{Bibliography}
\bibliographystyle{plain}
\bibliography{literature}

\end{document}